\documentclass[reqno]{amsart}
\usepackage{setspace}
\usepackage{graphicx}
\usepackage{amsmath}
\usepackage{amsfonts}
\usepackage{amssymb}

\begin{document}

\centerline{\Large \bf {\huge \bf C}LASSIFICATION OF THE SUBLATTICES OF A LATTICE}

\bigskip\medskip
\centerline{{\Large  C}HUANMING {\Large Z}ONG}

\vspace{0.8cm}
\centerline{\begin{minipage}{12.8cm}
{\bf Abstract.} In 1945-46, C. L. Siegel proved that an $n$-dimensional lattice $\Lambda $ of determinant ${\rm det}(\Lambda )$ has at most $m^{n^2}$ different sublattices of determinant $m\cdot {\rm det}(\Lambda )$. In 1997, the exact number of the different sublattices of index $m$ was determined by Baake. This paper presents a systematic treatment for counting the sublattices and deduces a formula for the number of the sublattice classes of determinant $m\cdot {\rm det}(\Lambda )$.
\end{minipage}}


\vspace{0.8cm}
\noindent
{\Large\bf 1. Introduction}

\bigskip\noindent
Let $\mathbb{Z}$ denote the set of all integers and let $\mathbb{E}^n$ denote the $n$-dimensional Euclidean space. If ${\bf a}_1$, ${\bf a}_2$,
$\ldots $, ${\bf a}_n$ are $n$ independent vectors in $\mathbb{E}^n$, then the discrete set
$$\Lambda =\left\{ \sum z_i{\bf a}_i:\ z_i\in \mathbb{Z}\right\}$$
is called an {\it $n$-dimensional lattice} generated by a {\it basis} $\{ {\bf a}_1, {\bf a}_2, \ldots , {\bf a}_n\}$. Assume that ${\bf a}_i=(a_{i1},$ $a_{i2},$ $\ldots,$ $a_{in})$, then the absolute value of the determinant of
$$A=\left(
\begin{array}{cccc}
a_{11}&a_{12}&\ldots &a_{1n}\\
a_{21}&a_{22}&\ldots &a_{2n}\\
\vdots &\vdots &\vdots &\vdots \\
a_{n1}&a_{n2}&\ldots &a_{nn}
\end{array}
\right)$$
is called the {\it determinant} of $\Lambda $. Usually, it is written as ${\rm det} (\Lambda )$. In fact, we also have
$${\rm det}(\Lambda )={\rm vol} (P),$$
where $P$ is the parallelopiped defined by
$$P=\left\{\sum \lambda_i {\bf a}_i:\ 0\le \lambda_i\le 1\right\}.$$

A subset $\Lambda^*$ of $\Lambda $ is called its {\it sublattice} if itself is an $n$-dimensional lattice as well. If $\{{\bf b}_1, {\bf b}_2, \ldots ,$ ${\bf b}_n \}$ is a basis of $\Lambda^*$, where ${\bf b_i}=(b_{i1}, b_{i2}, \ldots , b_{in})$, then we have
$${\bf b}_i=d_{i1}{\bf a}_1+d_{i2}{\bf a}_2+\ldots +d_{in}{\bf a}_n, \hspace{0.6cm} d_{ij}\in \mathbb{Z}.$$
Let $B$ denote the $n\times n$ matrix with elements $b_{ij}$ and let $D$ denote the $n\times n$ matrix with elements $d_{ij}$. Then we get
$$B=DA$$
and therefore
$${\rm det}(\Lambda^*)=m\cdot {\rm det}(\Lambda ),$$
where $m$ is the absolute value of the determinant of $D$. Usually, we call $m$ the {\it index} of $\Lambda^*$ in $\Lambda $.

The structures and representations of the sublattices have been studied by many authors such as Minkowski, Siegel, Cassels, Hlawka, Rogers and Schmidt. Many results and their applications can be found in classic references such as \cite{cassels,gruber-1,gruber-2,minkowski,siegel}. Particular sublattices have been studied by \cite{bernstein,fukshansky-1,fukshansky-2,schmidt-1,schmidt-2}.

\medskip
Let $\Lambda $ be an $n$-dimensional lattice, let $m$ be a positive integer, let $f_n(m)$ denote the number of the different sublattices of $\Lambda $ with index $m$, and let $f^*_n(m)$ denote the number of the different sublattice classes of $\Lambda $ with index $m$.

In 1945-46, C. L. Siegel gave a series of lectures on Geometry of Numbers at New York University. His lecture notes \cite{siegel} contained the first upper bound for
$f_n(m)$, namely
$$f_n(m)\le m^{n^2}.\eqno(1)$$
Since the lecture notes was published only in 1989, this result and many others were neglected.
In 1959, J. W. S. Cassels \cite{cassels} presented some basic result about the structures of the bases of the sublattices. In 1997, M. Baake \cite{baake} deduced the following formula based on a recursion in Algebra
$$f_n(m)=\sum_{d_1d_2\ldots d_n=m}d_1^0d_2^1\ldots d_n^{n-1}.\eqno(2)$$
Clearly, both Cassels and Baake were unaware of Siegel's work. Assume that
$$m=p_1^{\alpha_1}\ldots p_\ell^{\alpha_\ell},\eqno(3)$$
where $p_i$ are prime numbers. Baake's formula was simplified by B. Gruber \cite{gruber} as
$$f_n(m)=\prod_{i=1}^\ell \prod_{j=1}^{\alpha_\ell} {{p_i^{j+n-1}-1}\over {p_i^j-1}}=\prod_{i=1}^\ell \prod_{j=1}^{n-1} {{p_i^{j+\alpha_i}-1}\over {p_i^j-1}}.\eqno(4)$$
In particular, when $p$ is a prime, it is interesting to notice that
$$f_n(p)=1+p+\ldots +p^{n-1}$$
and
$$f_2(p^\ell)=1+p+\ldots +p^\ell .$$

\medskip
Let $k$ be a positive integer and let $p_n(k)$ denote the partition number of $k$ into $n$ parts. In other words, $p_n(k)$ is the number of the integer solutions for
$$\left\{
\begin{array}{ll}
x_1+x_2+\ldots +x_n=k,&\\
x_1\ge x_2\ge \ldots \ge x_n\ge 0.&
\end{array}\right.$$
The purpose of this paper is to present a systematic treatment on this topic, to complete both the statement and the proof. First, we present detailed proofs for (2) and (4). Then, we prove the following classification theorem.

\medskip\noindent
{\bf Theorem Z.} {\it If $m=p_1^{\alpha_1}\ldots p_\ell^{\alpha_\ell},$ where $p_i$ are prime numbers, then we have}
$$f_n^*(m)= \prod_{i=1}^\ell p_n(\alpha_i).$$

\smallskip
\noindent
{\bf Remark 1.} {\it When $m=p_1p_2\ldots p_\ell$, where $p_1$, $p_2$, $\ldots $, $p_\ell $ are pairwise distinct primes, we have
$$f_n(m)=\prod_{i=1}^\ell \sum_{j=0}^{n-1}p_i^j$$
and
$$f_n^*(m)=1.$$
Then, all the sublattices of index $m$ are equivalent to each others under unimodular transformations.}

\vspace{0.8cm}
\noindent
{\Large\bf 2. C. L. Siegel's Upper Bound}

\bigskip\noindent
Siegel's upper bound (1) was obtained in 1945-46. However, it was published only in 1989 in his lecture notes by Chandrasekharan \cite{siegel}.
So, this beautiful result has been neglected by almost all authors on related topics. For this reason, we reproduce it here. First of all, let us introduce a well-known basic lemma which can be found in every book on lattices.

\medskip\noindent
{\bf Lemma 1.} {\it Let $\{ {\bf a}_1, {\bf a}_2, \ldots , {\bf a}_n\}$ be a basis of an $n$-dimensional lattice $\Lambda$. Assume that ${\bf u}_1$, ${\bf u}_2$, $\ldots $, ${\bf u}_n$ are $n$ linear independent vectors in $\mathbb{E}^n$ with
$$ {\bf u}_i=u_{i1}{\bf a}_1+u_{i2}{\bf a}_2+\ldots +u_{in}{\bf a}_n,\hspace{0.5cm} i=1, 2, \ldots , n.$$
Then, $\{ {\bf u}_1, {\bf u}_2, \ldots , {\bf u}_n\}$ is also a basis of $\Lambda $ if and only if $U=\left( u_{ij}\right)$ is an $n\times n$ unimodular matrix.}

\medskip\noindent
{\bf Theorem 1 (Siegel \cite{siegel}).} {\it Assume that $\Lambda$ is an $n$-dimensional lattice and $m$ is a positive integer. Then $\Lambda $ has at most $m^{n^2}$ different sublattices of index $m$. In other words, we have}
$$f_n(m)\le m^{n^2}.$$

\noindent
{\bf Proof.} Assume that $\{ {\bf a}_1, {\bf a}_2, \ldots , {\bf a}_n\}$ is a basis of $\Lambda $. If $\Lambda^*$ is a sublattice of $\Lambda $ of index $m$ with a basis $\{ {\bf u}_1, {\bf u}_2, \ldots , {\bf u}_n\}$, then we have
$${\bf u}_i=u_{i1}{\bf a}_1+u_{i2}{\bf a}_2+\ldots +u_{in}{\bf a}_n, \hspace{0.4cm}i=1, 2, \ldots, n,\eqno(5)$$
where all $u_{ij}$ are integers and ${\rm det}(u_{ij})=\pm m$. For convenience, we denote the $n\times n$ matrix $\left(u_{ij}\right)$ by $U$.
If $\Lambda^\bullet$ is another sublattice of $\Lambda $ of index $m$ with a basis $\{ {\bf v}_1, {\bf v}_2, \ldots , {\bf v}_n\}$, then we have
$${\bf v}_i=v_{i1}{\bf a}_1+v_{i2}{\bf a}_2+\ldots +v_{in}{\bf a}_n, \hspace{0.4cm}i=1, 2, \ldots, n,\eqno(6)$$
where all $v_{ij}$ are integers and ${\rm det}(v_{ij})=\pm m$. We denote the $n\times n$ matrix $\left(v_{ij}\right)$ by $V$.

Clearly, it follows by (5) and (6) that the matrix that transforms $\{ {\bf v}_1, {\bf v}_2, \ldots , {\bf v}_n\}$ into $\{ {\bf u}_1, {\bf u}_2, \ldots ,$ ${\bf u}_n\}$ is $UV^{-1}$. In other words, if $W=UV^{-1}=\left(w_{ij}\right)$, we have
$${\bf u}_i=w_{i1}{\bf v}_1+w_{i2}{\bf v}_2+\ldots +w_{in}{\bf v}_n, \hspace{0.4cm}i=1, 2, \ldots, n.\eqno(7)$$

Now, we proceed to show that if
$$u_{ij}\equiv v_{ij} \quad ({\rm mod}\ m)$$
hold for all $i, j=1, 2, \ldots , n,$ then $\Lambda^*$ is identical with $\Lambda^\bullet$. Clearly $mV^{-1}$ is an integer matrix. Then, since $U\equiv V\ ({\rm mod}\ m)$, we have
$$mW=mUV^{-1}\equiv  mVV^{-1}\equiv m E \equiv O \quad ({\rm mod}\ m),\eqno(8)$$
where $E$ is the $n\times n$ unit matrix and $O$ is the $n\times n$ zero matrix. This means that all elements of $mW$ are divisible by $m$ and
therefore all elements of $W$ are integers. On the other hand, we have
$${\rm det}(W)={\rm det}(UV^{-1})=\pm {m\over m}=\pm 1.\eqno(9)$$
Thus, $W$ must be a unimodular matrix. Then it follows by Lemma 1 that $\Lambda^*$ is identical with $\Lambda^\bullet$.

This shows that there are at most $m$ possible values for any element of $U$, such that the corresponding sublattices of $\Lambda $ are different.
Since $U$ has $n^2$ elements, the total number of possibilities for $U$ is $m^{n^2}$. In other words,
$$f_n(m)\le m^{n^2}.$$
The theorem is proved. \hfill{$\square$}

\vspace{0.8cm}
\noindent
{\Large\bf 3. The Sublattices of Given Index}

\bigskip\noindent
In 1907, Minkowski \cite{minkowski} studied the relation between the bases of a three-dimensional lattice and its sublattices. Afterwards, his result was generalized into arbitrary dimensions (see \cite{cassels} or \cite{gruber-1}) as following. Assume that $\Lambda^*$ is a sublattice of an $n$-dimensional $\Lambda $. If $\{ {\bf u}_1, {\bf u}_2, \ldots , {\bf u}_n\}$ is a basis of $\Lambda^*$, then $\Lambda $ has a basis $\{ {\bf a}_1, {\bf a}_2, \ldots , {\bf a}_n\}$ satisfying
$${\bf u}_i=u_{i1}{\bf a}_1+u_{i2}{\bf a}_2+\ldots +u_{ii}{\bf a}_i,\quad i=1, 2, \ldots , n,\eqno(10)$$
where $u_{ii}>0$ and $0\le u_{ij}<u_{ii}$ for all $j<i$.

It is rather unexpected that the following inverse of this result is also true. It can be found in both \cite{cassels} and \cite{gruber-1}, neither of them indicated further reference.

\medskip\noindent
{\bf Lemma 2 (Cassels \cite{cassels}).} {\it Assume that $\Lambda $ is an $n$-dimensional lattice with a basis $\{ {\bf a}_1, {\bf a}_2, \ldots , {\bf a}_n\}$. If $\Lambda^*$ is a sublattice of $\Lambda $ of index $m$, then $\Lambda^*$ has a basis $\{ {\bf u}_1, {\bf u}_2, \ldots , {\bf u}_n\}$ satisfying
$${\bf u}_i=u_{i1}{\bf a}_1+u_{i2}{\bf a}_2+\ldots +u_{ii}{\bf a}_i,\quad i=1, 2, \ldots , n$$
and
$$m=u_{11}u_{22}\ldots u_{nn},$$
where $u_{ii}>0$ and $0\le u_{ij}<u_{jj}$ for all $j<i$.}

\medskip
Clearly, this lemma provides a mean to count the number of the different sublattices of given index $m$. To do the explicit counting, we need another simple result.

\medskip\noindent
{\bf Lemma 3.} {\it Assume that $\Lambda $ is an $n$-dimensional lattice with a basis $\{ {\bf a}_1, {\bf a}_2, \ldots , {\bf a}_n\}$ and $m$ is a positive integer. Let ${\bf u}_1, {\bf u}_2, \ldots , {\bf u}_n $ be $n$ linearly independent vectors satisfying
$${\bf u}_i=u_{i1}{\bf a}_1+u_{i2}{\bf a}_2+\ldots +u_{ii}{\bf a}_i,\quad i=1, 2, \ldots , n$$
and
$$m=u_{11}u_{22}\ldots u_{nn},$$
where all $u_{ij}$ are integers, $u_{ii}>0$ and $0\le u_{ij}<u_{jj}$ for all $j<i$, let ${\bf v}_1, {\bf v}_2, \ldots , {\bf v}_n $ be $n$ linearly independent vectors satisfying
$${\bf v}_i=v_{i1}{\bf a}_1+v_{i2}{\bf a}_2+\ldots +v_{ii}{\bf a}_i,\quad i=1, 2, \ldots , n$$
and
$$m=v_{11}v_{22}\ldots v_{nn},$$
where all $v_{ij}$ are integers, $v_{ii}>0$ and $0\le v_{ij}<v_{jj}$ for all $j<i$, let $\Lambda^*$ be the sublattice of $\Lambda $ generated by $\{ {\bf u}_1, {\bf u}_2, \ldots , {\bf u}_n \}$, and let $\Lambda^\bullet$ be the sublattice of $\Lambda $ generated by $\{ {\bf v}_1, {\bf v}_2, \ldots , {\bf v}_n \}$. Then the two sublattices $\Lambda^*$ and $\Lambda^\bullet$ are identical if and only if}
$$u_{ij}=v_{ij},\quad 1\le j\le i \le n.$$

\smallskip\noindent
{\bf Proof.} The if part is obvious. Now, let us prove the only if part.

Let $U$ denote the $n\times n$ matrix with elements $u_{ij}$, $i,j=1, 2, \ldots , n$, where $u_{ij}=0$ for all $j>i$,
let $V$ denote the $n\times n$ matrix with elements $v_{ij}$, $i,j=1, 2, \ldots , n$, where $v_{ij}=0$ for all $j>i$, and define
$$W=UV^{-1}=\left(w_{ij}\right).\eqno(11)$$
It is easy to see that $\Lambda^*=\Lambda^\bullet$ if and only if $W$ is a unimodular matrix.

By (11) we have
$$WV=U.\eqno(12)$$
Then, by comparing both sides of (12) for $u_{1n}$, $u_{1,n-1}$, $\ldots $, $u_{11}$, we get
$$\left\{
\begin{array}{ll}
w_{11}v_{1n}+w_{12}v_{2n}+\ldots +w_{1n}v_{nn}=0, &\\
w_{11}v_{1,n-1}+w_{12}v_{2,n-1}+\ldots +w_{1n}v_{n,n-1}=0, &\\
\hspace{2.5cm}\ldots ,&\\
w_{11}v_{11}+w_{12}v_{21}+\ldots +w_{1n}v_{n1}=u_{11}&
\end{array}\right.$$
and thus
$$\left\{
\begin{array}{ll}
w_{1n}=w_{1,n-1}=\ldots =w_{12}=0,&\\
w_{11}v_{11}=u_{11}.&
\end{array}\right.\eqno(13)$$
Repeating this process for $u_{2i},$ $u_{3i}$, $\ldots,$ $u_{ni}$ successively, we get
$$\left\{
\begin{array}{ll}
w_{ij}=0,\quad i<j\le n, &\\
w_{ii}v_{ii}=u_{ii}, \quad i=1, 2, \ldots , n.&
\end{array}\right.\eqno(14)$$

If $W$ is a unimodular matrix, all its elements are integers, it follows by (14) and the assumption
$$m=u_{11}u_{22}\ldots u_{nn}= v_{11}v_{22}\ldots v_{nn}$$
that
$$w_{11}=w_{22}=\ldots =w_{nn}=1.\eqno(15)$$
Then, by comparing both sides of (12) for $u_{21}$, $u_{32}$, $\ldots ,$ $u_{n,n-1}$, we get
$$w_{i+1,i}v_{ii}+v_{i+1,i}=u_{i+1,i},\qquad i=1, 2, \ldots, n-1.\eqno(16)$$
If $w_{i+1, i}\not= 0,$ by (16) we get
$$w_{i+1,i}v_{ii}=u_{i+1,i}-v_{i+1,i},$$
which contradicts the assumptions that $0\le u_{i+1, i}<u_{ii}=v_{ii}$ and $0\le v_{i+1, i}<v_{ii}$. Thus, we must have
$$\left\{
\begin{array}{ll}
w_{i+1,i}=0, &\\
u_{i+1,i}=v_{i+1,i}&
\end{array}\right.\eqno(17)$$
for all $i=1,2, \ldots, n-1.$

Inductively, assume that
$$w_{i+j,i}=0\eqno(18)$$
holds for all $1\le j\le k-1<n-1$ and $i=1, 2, \ldots, n-j,$ by comparing both sides of (12) for $u_{i+k,i}$, $i=1, 2, \ldots, n-k,$ similar to (16)
we can get
$$w_{i+k,i}=0,\quad i=1,2, \ldots, n-k.\eqno(19)$$

As a conclusion, we obtain that, if $W$ is a unimodular matrix, it must be the $n\times n$ unit matrix. In other words, if $\Lambda^*=\Lambda^\bullet$, then $U=V.$ The theorem is proved. \hfill{$\square$}

\medskip
Clearly, an $n$-dimensional lattice is a free module of rank $n$ over $\mathbb{Z}$. By studying the algebraic structures of the submodules it was shown (see \cite{scheja}) that
$$f_n(m)=\sum_{d\mid m}d\cdot f_{n-1}(d).\eqno(20)$$
In 1997, it was deduced from (20) by Baake \cite{baake} that
$$f_n(m)=\sum_{d_1d_2\ldots d_n=m}d_1^0d_2^1\ldots d_n^{n-1}.\eqno(21)$$

In fact, Baake's formula can be easily deduced from Lemma 2 and Lemma 3. Gruber \cite{gruber} did realize this possible connection and simplified (21). However, he neglected the necessity of Lemma 3.

\medskip\noindent
{\bf Theorem 2 (Baake \cite{baake}, Gruber \cite{gruber}).} {\it If $m=p_1^{\alpha_1}p_2^{\alpha_2}\ldots p_\ell^{\alpha_\ell}$, where $p_i$ are prime numbers and $\alpha_i$ are positive integers, then}
$$f_n(m)=\sum_{d_1d_2\ldots d_n=m}d_1^0d_2^1\ldots d_n^{n-1}=\prod_{i=1}^\ell \prod_{j=1}^{\alpha_i}{{p_i^{j+n-1}-1}\over {p_i^j-1}}
=\prod_{i=1}^\ell \prod_{j=1}^{n-1}{{p_i^{j+\alpha_i}-1}\over {p_i^j-1}}.$$

\medskip\noindent
{\bf Remark 2.} Noticing that
$$\bigl(p_i^{j+n-1}-1\bigr)/\bigl(p_i^j-1\bigr)\le p_i^n$$
and
$$\bigl(p_i^{j+\alpha_i}-1\bigr)/\bigl(p_i^j-1\bigr)\ge p_i^{\alpha_i},$$
one can easily deduce that
$$m^{n-1}\le f_n(m)\le m^n.$$
Comparing with Theorem 1, it is interesting to see that Siegel's upper bound is far away from the exact values of $f_n(m)$.

\vspace{0.8cm}
\noindent
{\Large\bf 4. Classification of the Sublattices of Given Index}

\bigskip\noindent
Let $\Lambda $ be an $n$-dimensional lattice in $\mathbb{E}^n$, and let $\Lambda^*$ and $\Lambda^\bullet$ be two sublattices of $\Lambda $. We say
that $\Lambda^*$ and $\Lambda^\bullet$ are {\it equivalent} if there is a linear transformation $\sigma $ satisfying both
$$\sigma (\Lambda )=\Lambda \eqno(22)$$
and
$$\sigma (\Lambda^* )=\Lambda^\bullet. \eqno(23)$$
Then, for convenience, we write $\Lambda^* \sim \Lambda^\bullet.$  Clearly, a linear transformation satisfying $\sigma (\Lambda )=\Lambda $ if and only if $\sigma $ is corresponding to a unimodular matrix.

\medskip\noindent
{\bf Example 1.} Let $\Lambda =\mathbb{Z}^2$ with ${\bf e}_1=(1,0)$ and ${\bf e}_2=(0,1)$, let $\Lambda^*$ be the sublattice generated by ${\bf u}_1={\bf e}_1$ and ${\bf u}_2=2{\bf e}_2$, and let $\Lambda^\bullet$ be the sublattice generated by ${\bf u}_1=2{\bf e}_1$ and ${\bf u}_2={\bf e}_2$. It is obvious that $\Lambda^*\not=\Lambda^\bullet$. Let $\sigma $ denote the linear transformation determined by $\sigma ({\bf e}_1)={\bf e}_2$ and $\sigma ({\bf e}_2)={\bf e}_1$, it can be verified that $\sigma (\Lambda )=\Lambda $ and $\sigma (\Lambda^*)=\Lambda^\bullet.$ Thus, we have $\Lambda^*\sim\Lambda^\bullet.$

\medskip
It is shown in Gruber \cite{gruber-1} that, if $\Lambda^*$ is a sublattice of $\Lambda $, then $\Lambda $ has a basis $\{ {\bf a}_1, {\bf a}_2, \ldots , {\bf a}_n\}$ and $\Lambda^*$ has a basis $\{ {\bf u}_1, {\bf u}_2, \ldots , {\bf u}_n\}$ such that
$${\bf u}_i=u_{ii}{\bf a}_i,\quad i=1, 2, \ldots , n,\eqno(24)$$
where $u_{ii}$ are suitable positive integers.

On page 26 of \cite{martinet}, Martinet wrote \lq\lq Let $M$ be an $R$-module and let $M'$ be a submodule of $M$, both having the same rank $n$.
(When $R=\mathbb{Z}$, this amounts to saying that $[ M:M']<\infty$.) There then exists a basis $B=\{ {\bf e}_1, {\bf e}_2, \ldots , {\bf e}_n\}$ for $M$ and nonzero elements $a_1$, $a_2$, $\ldots $, $a_n$ of $R$ such that $B'=\{a_1{\bf e}_1, a_2{\bf e}_2, \ldots , a_n{\bf e}_n\}$ is a basis for $M'$, and $a_i$ divides $a_{i-1}$ for $2\le i\le n$." This implies that $u_{ii}$ divides $u_{i-1,i-1}$ in (24).

For the completeness, we restate this result as Lemma 4 in the following and give a detailed proof.

\medskip\noindent
{\bf Lemma 4.} {\it If $\Lambda^*$ is a sublattice of $\Lambda $, then $\Lambda $ has a basis $\{ {\bf a}_1, {\bf a}_2, \ldots , {\bf a}_n\}$ and $\Lambda^*$ has a basis $\{ {\bf u}_1, {\bf u}_2, \ldots , {\bf u}_n\}$ such that
$${\bf u}_i=u_{ii}{\bf a}_i,\quad i=1, 2, \ldots , n,$$
where all $u_{ii}$ are positive integers satisfying $u_{ii}\mid u_{i-1,i-1}$ for all $2\le i\le n$. }

\medskip\noindent
{\bf Proof.} Assume that $\{{\bf e}_1, {\bf e}_2, \ldots , {\bf e}_n\}$ is a basis for $\Lambda $ and $\{{\bf v}_1, {\bf v}_2, \ldots , {\bf v}_n\}$ is a basis for $\Lambda^*$. Then, we have
$${\bf v}_i=v_{i1}{\bf e}_1+v_{i2}{\bf e}_2+\ldots +v_{in}{\bf e}_n, \quad i=1, 2, \ldots , n. \eqno(25)$$
For convenience, let $\overline{\bf X}$ denote the $n\times 1$ matrix with elements ${\bf x}_1$, ${\bf x}_2$, $\ldots $, ${\bf x}_n$ and let $X$ denote the $n\times n$ matrix with elements $x_{ij}$. Then, one can rewrite (25) as
$$\overline{\bf V}=V\overline{\bf E}.\eqno(26)$$

If $\{ {\bf u}_1, {\bf u}_2, \ldots , {\bf u}_n\}$ is another basis for $\Lambda^*$ such that
$$\overline{\bf V}=U_1\overline{\bf U},\eqno(27)$$
where $U_1$ is an $n\times n$ unimodular matrix, and $\{ {\bf a}_1, {\bf a}_2, \ldots , {\bf a}_n\}$ is another basis for $\Lambda$ such that
$$\overline{\bf E}=U_2\overline{\bf A},\eqno(28)$$
where $U_2$ is an $n\times n$ unimodular matrix. Then, it follows by (26), (27) and (28) that
$$\overline{\bf U}=U_1^{-1}VU_2\overline{\bf A}.\eqno(29)$$

It is known in Algebra (see Chapter 14 of Hua \cite{hua}) that, for a given integer matrix $V$ there are two suitable unimodular matrices $U_1$ and $U_2$ such that
$$U_1^{-1}VU_2=\left(
\begin{array}{cccc}
u_{11}&0&\ldots &0\\
0&u_{22}&\ldots &0\\
\vdots &\vdots&\ddots &\vdots\\
0&0&\ldots&u_{nn}
\end{array}\right),$$
where $u_{ii}\mid u_{i-1,i-1}$ for all $2\le i\le n$. Then, by (29) we have
$${\bf u}_i=u_{ii}{\bf a}_i,\quad i=1, 2, \ldots, n.$$
The lemma is proved. \hfill{$\square$}

\medskip\noindent
{\bf Lemma 5.} {\it  Assume that $\Lambda^*$ and $\Lambda^\bullet$ are two sublattices of an $n$-dimensional lattice $\Lambda $. If $\{ {\bf u}_1, {\bf u}_2, \ldots ,$ ${\bf u}_n\}$ is a basis of $\Lambda^*$ and $\{ {\bf a}_1, {\bf a}_2, \ldots , {\bf a}_n\}$ is a basis of $\Lambda$ such that
$${\bf u}_i=u_{ii}{\bf a}_i,\quad i=1, 2, \ldots , n,$$
where $u_{ii}$ are positive integers satisfying $u_{ii}\mid u_{i-1,i-1}$ for all $2\le i\le n,$  and $\{ {\bf v}_1, {\bf v}_2, \ldots , {\bf v}_n\}$ is a basis of $\Lambda^\bullet$ and $\{ {\bf b}_1, {\bf b}_2, \ldots , {\bf b}_n\}$ is a basis of $\Lambda$ such that
$${\bf v}_i=v_{ii}{\bf b}_i,\quad i=1, 2, \ldots , n,$$
where $v_{ii}$ are positive integers satisfying $v_{ii}\mid v_{i-1,i-1}$ for all $2\le i\le n.$ Then, $\Lambda^*\sim \Lambda^\bullet$ if and only if}
$$u_{ii}=v_{ii},\quad i=1, 2, \ldots, n.$$

\medskip\noindent
{\bf Proof.} If $u_{ii}=v_{ii}$ hold for all $i=1, 2, \ldots, n.$ Let $\sigma $ be the linear transformation defined by
$$\sigma ({\bf a}_i)={\bf b}_i,\quad i=1, 2, \ldots , n,$$
then we have
$$\sigma ({\bf u}_i)=\sigma (u_{ii}{\bf a}_i)=u_{ii}{\bf b}_i={\bf v}_i$$
for all $i=1, 2, \ldots , n$ and thus
$$\sigma (\Lambda^*)=\Lambda^\bullet .$$

On the other hand, if $\Lambda^*\sim \Lambda^\bullet $ with a suitable $\sigma$, then we have
$$\overline{\bf U}=U\overline{\bf A},\eqno(30)$$
$$\overline{\bf V}=V\overline{\bf B},\eqno(31)$$
$$\sigma \left(\overline{\bf U}\right)=W\overline{\bf V}\eqno(32)$$
and
$$\sigma \left(\overline{\bf A}\right)=T\overline{\bf B},\eqno(33)$$
where $u_{ij}=0$ for all $i\not= j$, $v_{ij}=0$ for all $i\not= j$, both $W$ and $T$ are suitable unimodular matrices.

It follows by $\sigma (\Lambda^*)=\Lambda^\bullet $, (30), (31), (32) and (33) that
$$\sigma \left(\overline{\bf U}\right)=U \sigma \left(\overline{\bf A}\right),$$
$$W\overline{\bf V} =UT\overline{B},$$
$$\overline{\bf V}=W^{-1}UT\overline{\bf B}=V\overline{\bf B}$$
and thus
$$V=W^{-1}UT.\eqno(34)$$
It is known in Algebra (see Chapter 14 of Hua \cite{hua}) that (34) implies
$$V=U.$$
Lemma 5 is proved. \hfill{$\square$}

\medskip\noindent
{\bf Proof of Theorem Z.}  Recall that $$m=p_1^{\alpha_1}p_2^{\alpha_2}\ldots p_\ell^{\alpha_\ell},$$
where $p_i$ are prime numbers. It follows by Lemma 4 and Lemma 5 that $f_n^*(m)$ is the number of the factorizations
$$m=d_1d_2\ldots d_n\eqno(35)$$
satisfying $d_j\mid d_{j-1}$ for all $2\le j\le n$. If
$$d_j=p_1^{\beta_{1j}}p_2^{\beta_{2j}}\ldots p_\ell^{\beta_{\ell j}},\eqno(36)$$
then we have
$$\left\{
\begin{array}{ll}
\sum_{j=1}^n\beta_{ij}=\alpha_i,&\\
\beta_{i1}\ge \beta_{i2}\ge \ldots \ge \beta_{in}\ge 0&
\end{array}\right.\eqno(37)$$
for all $i=1, 2, \ldots, \ell .$ Clearly (37) has $p_n(\alpha_i)$ solutions and each solution corresponds to one factorization of (35). Thus, we have
$$f_n^*(m)=\prod_{i=1}^\ell p_n(\alpha_i).$$
Theorem Z is proved. \hfill{$\square$}

\medskip\noindent
{\bf Remark 3.} The partition function $p_n(k)$ has been studied by many authors. See Andrews and Eriksson \cite{andrews} for references.

\vspace{0.8cm}\noindent
{\bf Acknowledgements.} For helpful correspondence, the author is grateful to Prof. Baake. This work is supported by 973 Program 2013CB834201.

\medskip
\bibliographystyle{amsplain}

\begin{thebibliography}{99}
\bibitem{andrews}G. E. Andrews and K. Eriksson, {\it Integer Partitions}, Cambridge University Press, 2004.
\bibitem{baake}M. Baake, Solution of the coincidence problem in dimensions $d\le 4$, {\it The Mathematics of Long-Range Aperiodic Order}, edited by R. V. Moody, 1997, 9-44. Dordrecht, Kluwer Academic Publishers. arXiv:math/0605222.
\bibitem{bernstein}M. Bernstein, N. J. A. Sloane, P. E. Wright, On sublattices of the hexagonal lattice, {\it Discrete Math.} {\bf 170} (1997), 29-39.
\bibitem{cassels}J. W. S. Cassels, {\it An Introduction to the Geometry of Numbers}, Springer-Verlag, Berlin, 1959.
\bibitem{fukshansky-1}L. Fukshansky, On distribution of well-rounded sublattices of $\mathbb{Z}^2$, {\it J. Number Theory}, {\bf 128} (2008), 2359-2393.
\bibitem{fukshansky-2}L. Fukshansky, On similarity classes of well-rounded sublattices of $\mathbb{Z}^2$, {\it J. Number Theory}, {\bf 129} (2009), 2530-2556.
\bibitem{gruber}B. Gruber, Alternative formulae for the number of sublattices, {\it Acta Cryst.} {\bf A53} (1997), 807-808.
\bibitem{gruber-1}P. M. Gruber, {\it Convex and Discrete Geometry}, Springer-Verlag, New York, 2007.
\bibitem{gruber-2}P. M. Gruber and C. G. Lekkerkerker, {\it Geometry of Numbers}, North-Holland, Amsterdam, 1987.
\bibitem{hua}L. K. Hua, {\it Introduction to Number Theory}, Springer, 1987.
\bibitem{martinet}J. Martinet, {\it Perfect Lattices in Euclidean Spaces}, Springer-Verlag, Berlin, 2003.
\bibitem{minkowski}H. Minkowski, {\it Diophantische Approximationen}, Teubner, Leipzig, 1907; Chelsea, New York, 1957.
\bibitem{scheja}G. Scheja and U. Storch, {\it Lehrbuch der Algebra}, Teil 2, Teubner, Stuttgart, 1988.
\bibitem{schmidt-1}W. M. Schmidt, The distribution of the sublattices of $\mathbb{Z}^n$, {\it Monatsh. Math.} {\bf 125} (1998), 37-81.
\bibitem{schmidt-2}W. M. Schmidt, Integer matrices, sublattices of $\mathbb{Z}^n$, and Frobenius numbers, {\it Monatsh. Math.} {\bf 178} (2015), 405-451.
\bibitem{siegel}C. L. Siegel, {\it Lectures on the Geometry of Numbers}, Springer-Verlag, Berlin, 1989.
\end{thebibliography}

\vspace{0.8cm}
\noindent
Chuanming Zong, Center for Applied Mathematics, Tianjin University, Tianjin 300072, China

\noindent
Email: cmzong@math.pku.edu.cn

\end{document}